\input amstex
\pageno=1

\magnification=1200
\loadmsbm
\loadmsam
\loadeufm
\input amssym
\UseAMSsymbols

\TagsOnRight

\hsize163 true mm
\vsize220 true mm
\voffset=20 true mm
\hoffset=0 true mm
\baselineskip 8.5 true mm plus0.4 true mm minus0.4 true  mm

\def\P{\partial}

\def\A{\alpha}
\def\B{\beta}

\def\X{\xi}
\def\F{\psi}
\def\SI{\sigma}

\vskip 0.4cm
\noindent
{\bf Witten solution of the Gelfand-Dikii hierarchy}

\noindent
{\bf S.M.Natanzon}

\vskip 0.3cm
\noindent
{\it Moscow State University, e-mail: natanzon\@mccme.ru}
\vskip 0.4cm

\noindent
{\bf Abstract.} We produce formulas, permiting to
find the coefficients of Taylor--series expanded of some im\-portant
solution of the Gelfand-Dikii hierarchy. By the Witten conjecture
these co\-ef\-ficients are numbers of intersection of
Mumford-Morita-Muller stable cohomolo\-gical classes of moduli
space of $n$-spin bundles on Riemann surfaces with punc\-tures.



\vskip 0.5cm
\noindent
{\bf 1. Introduction}
\vskip 0.4cm

\noindent
The $n$-Gelfand-Dikii hierarchy ($n$-K$\P$V
hierarchy) usually is described in term of formal pseudo-differential
operators. Consider functions $u_j=u_j(x)\ (j=0,\dots, n-2)$
of infinite number of variable $x=(x_1, x_2,\dots .)$. Let us
$\P_i=\frac{\P}{\P x_i},\quad \P=\P_1$ and $L=\P^n+u_{n-2}\P^{n-2}+
\dotsb+u_1\P+u_0$. Denote by $L^{\frac in}$ the pseudo-differential
operator such that $(L^{\frac in})^n= L^i$. Let $[L^{\frac in}]_+$
be its differential part. Then the $n$-Gelfand-Dikii hierarchy
[GD] is the system of differential equations on $u_i$, which follow
from infinite system of equations
$$\P_i L=[[L^{\frac in}]_+,L]\quad (i=1,2,...).$$

Any solution $(u_0,\dots, u_{n-2})$ is generated by a function
$v(x_1, x_2,\dots)$. Among these solu\-tions there exists a remarkable
solution $W$, which satisfies the string equation,
generates a vacuum vector of $W$--algebra and
has a representation in a form of matrix integral [AM]. We call
it Witten's solution because according to the Witten conjecture
the function $F(x_1, x_2, x_3,\dots) = W(x_1,\frac{x_2}{2},
\frac{x_3}{3}, \dots)$ is the generating function
for numbers of intersec\-tion of
Mumford-Morita-Muller stable cohomolo\-gical classes of moduli
space of $n$-spin bundles on Riemann surfaces with punctures [W].
This conjecture was proved by Kontse\-vich for $n=2$ [Ko] and by Witten
himself for surfaces of genus 0 [W].

In this paper we find recurrence relation between coefficients of
Taylor series of W. This reduces the Witten conjecture to conjecture
that $n$-spin Mamford-Morita-Muller numbers satisfy to the same
relations. These relations give also an algorithm for calcula\-tion of
$n$-spin Mamford-Morita-Muller numbers in assuming that the Witten
conjecture is true. Moreover we prove that $F(x_1, x_2, \dots,
x_{n-1},0,0,\dots) =W(x_1,\frac{x_2}{2},\dots, \frac{x_n}{n-1},0,0,\dots)$
and thus the Witten conjecture is true for numbers $(c_D(\nu),
\bar{M}'_{g,s})$ in notation [W].

We prove also that the solution $W$ has a representation
$W=\sum\limits^\infty_{g=0} W^g$, where $W^g$ are quasihomohenious
series of degrees $(1-g)(2+\frac 2n)$.  This is some indirect
corroboration of the Witten conjecture because according to [W]
the function $F$ has a representation $F=\sum\limits^\infty_{g=0} F^g$,
where $F^g$ are quasihomohenious series of degrees
$(1-g)(2+\frac 2n)$, corresponding to surfaces of genius $g$.

We prove that $W^0(x_1, \dots, x_{n-1}, 0,\dots)=W(x_1, \dots, x_{n-1},
0, \dots)$ and therefore \linebreak
$F^0(x_1, \dots, x_{n-1}, 0, \dots)=F(x_1, \dots,
x_{n-1}, 0, \dots)$. Last function is polynomial
solution of WDVV equa\-tions. Some simple
formulas for calculation of this solution was found in [N2].

Organisation of the paper is as follow. In sect 2-4 we following by
[DN, N1] represent KP hierarchy as a system of differential equation
for $v=-\ln\tau$. In section 5 we prove that the $n$-hierarchy of
Gelfand-Dikii is equivalent of a system of differential equations
in a form
$$\P_{i_1}\dotsb\P_{i_k} v=\sum^\infty_{m=1} N^m_{i_1\dotsb i_k}
\pmatrix s_1 &\ldots & s_m\\t_1 &\ldots & t_m\endpmatrix
\P_{s_1}\P^{t_1}_1v\dotsb \P_{s_m}\P^{t_m}_1v, \tag 1.1$$
where $i_j, t_i\geqslant 1, s_i<n$ and $N^m_{i_1\dotsb i_k}
\pmatrix s_1 &\ldots & s_m\\t_1 &\ldots & t_m\endpmatrix$ are
rational constants.

In section 6 we investigate the Witten solution W of the system (1.1).

The author thanks B.A.Dubrovin for fruitful discussions.

\vskip 0.4cm
\noindent
{\bf 2. Combinatorial lemma}
\vskip 0.4cm

\noindent
For natural $s, i_1,...,i_n$ and integer not negative
$j_1,...,j_n$ defind $P_s\pmatrix i_1&\ldots&i_n\\j_1&\ldots&j_n
\endpmatrix$ by recu\-rence formulas:
$$1) P_s\pmatrix i_1&\ldots&i_n\\0&\ldots&0\endpmatrix=0; \quad
2) P_s\pmatrix i\\j\endpmatrix=C^j_s
\quad \text{for}\ j>0;$$
$$3) P_s\pmatrix i_1&\ldots&i_n\\j_1&\ldots&j_n\endpmatrix=\frac{1}{n!}
C_s^{j_1+\dotsb+j_n}\frac{(j_1+\dotsb+j_n)!}{j_1!\dotsb j_n!}-$$
$$\sum^{n-1}_{q=1}P_s
\pmatrix i_1 &\ldots&i_q\\j_1 &\ldots&j_q\endpmatrix\frac{1}{(n-q)!}
C^{j_{q+1}+\dotsb+j_n}_{s-(i_1+\dotsb+i_q+j_1+\dotsb+j_q)}
\frac{(j_{q+1}+\dotsb+j_n)!}{j_{q+1}!\dotsb j_n!}$$
for $(j_1,...,j_n)\ne (0,...,0)$.

Let $\bmatrix i_1 &\ldots&i_n\\j_1 &\ldots&j_n\endbmatrix$ be the set of
all matrices, which appear from $\pmatrix i_1 &\ldots&i_n\\j_1 &\ldots&j_n
\endpmatrix$ by permu\-tation of columns.
Let $\Vmatrix i_1 &\ldots&i_n\\j_1 &\ldots&j_n\endVmatrix$ be the number
of such matrices. Put us
$$P_s\bmatrix i_1 &\ldots&i_n\\j_1 &\ldots&j_n\endbmatrix=\sum
P_s\pmatrix a_1 &\ldots&a_n\\b_1 &\ldots&b_n\endpmatrix,$$ where the sum is
taken by all
$$\pmatrix a_1 &\ldots&a_n\\b_1 &\ldots&b_n\endpmatrix\in
\bmatrix i_1 &\ldots&i_n\\j_1 &\ldots&j_n\endbmatrix$$

\vskip 0.4cm
\noindent
{\bf Lemma 2.1.} {\sl Let $m>0, k>0$ and $j_n\geqslant 1$
for $n\leqslant m$. Then
$$P_s\bmatrix i_{1} & \ldots &i_m & i_{m+1} &\ldots &i_{m+k}\\
j_1 &\ldots &j_m & 0 &\ldots & 0\endbmatrix=$$
$$=\cases 0, \quad \text{if $s\geqslant i_1+...+i_m+j_1+...+j_m$,}\\
\frac{1}{k!}\Vmatrix i_{m+1}&\ldots& i_{m+k}\\0&\ldots &0\endVmatrix
P_s\bmatrix i_1&\ldots&i_m\\j_1&\ldots &j_m\endbmatrix, \quad
\text{if $s< i_1+...+i_m+j_1+...+j_m.$} \endcases$$}

\vskip 0.4cm
\noindent
Proof: Prove at first the lemma for $m=1$ using induction by k.
For $m=k=1$
$$P_s\bmatrix i_1 \ i_2\\j_1 \ 0\endbmatrix=
P_s\pmatrix i_1 \ i_2\\j_1 \ 0\endpmatrix+
P_s\pmatrix i_2 \ i_1\\0 \ j_1\endpmatrix=$$
$$=\frac 12 C^{j_1}_s-P_s\pmatrix i_1\\j_1\endpmatrix\cdot
C^0_{s-(i_1+j_1)}+\frac 12 C^{j_1}_s-P_s\pmatrix i_2\\0\endpmatrix
\cdot C^0_{s-i_2}=C^{j_1}_s-C^{j_1}_s C^0_{s-(i_1+j_1)}=$$
$$=\cases 0, &\text{if $s\geqslant i_1+j_1,$}\\
C^{j_1}_s=P_s\bmatrix i_1\\ j_1\endbmatrix, &\text{if
$s<i_1+j_1.$}\endcases$$
Prove the lemma for $m=1, k=N$, considering that it is proved for
$m=1, k<N$. If $s\geqslant i_1+j_1$, then
$$P_1\bmatrix i_1&i_2&\ldots & i_{k+1}\\j_1&0&\ldots & 0\endbmatrix=$$
$$=\frac{1}{k!}C^{j_1}_s
\Vmatrix i_2 &\ldots & i_{k+1}\\0&\ldots &  0\endVmatrix-
P_s\pmatrix i_1\\j_1\endpmatrix\cdot C^0_{s-(i_1+j_1)}\cdot
\frac{1}{k!}\Vmatrix i_2 &\ldots & i_{k+1}\\0 &\ldots & 0\endVmatrix=0.$$
If $s<i_1+j_1$ than
$$P_1\bmatrix i_1&i_2 &\ldots & i_{k+1}\\j_1 & 0 &\ldots &  0\endbmatrix=
\frac{1}{k!}C^{j_1}_s
\Vmatrix i_2 &\ldots & i_{k+1}\\0 &\ldots & 0\endVmatrix-
A \cdot C^0_{s-(i_1+j_1)}=$$
$$=\frac{1}{k!}\Vmatrix i_2 &\ldots & ,i_{k+1}\\0 &\ldots & 0\endVmatrix
P_s\bmatrix i_1\\ j_1\endbmatrix.$$ Thus the lemma is proved
for $m=1$.

Prove the lemma for $m=N$ considering that it is proved for $m<N$.
Then $$P_s\bmatrix i_1 &\ldots & i_m&i_{m+1} &\ldots & i_{m+k}\\j_1
 &\ldots & j_m &  0 &\ldots &  0\endbmatrix=$$
$$=\sum_{\pmatrix \A_1 &\ldots & \A_m\\ \B_1 &\ldots & \B_m\endpmatrix\in
\bmatrix i_1 &\ldots & i_m\\j_1 &\ldots & j_m\endbmatrix}\Bigl(\frac{1}
{(m+k)!} C^{\B_1+\dotsb+\B_m}_s\frac{(\B_1+\dotsb+\B_m)!}{\B_1!\dotsb \B_m!}
C^k_{m+k}\cdot$$
$$\cdot\Vmatrix i_{m+1} &\ldots & i_{m+k}\\ 0 &\ldots & 0\endVmatrix-$$
$$\sum_{q=1}^m P_s \pmatrix \A_1 &\ldots & \A_q\\
\B_1 &\ldots & \B_q\endpmatrix
\frac{1}{(m+k-q)!} C^{\B_{q+1}+\dotsb+\B_m}_{s-(\A_1+\dotsb
+\A_q+\B_1+\dotsb +\B_q)}
\frac{(\B_{q+1}+\dotsb+\B_m)!}{\B_{q+1}!\dotsb \B_m!}
C^k_{m+k-q}\cdot$$
$$\cdot\Vmatrix i_{m+1} &\ldots & i_{m+k}\\ 0 &\ldots & 0\endVmatrix\Bigr)=
P_s\bmatrix i_1 &\ldots & i_m\\ j_1 &\ldots & j_m\endbmatrix\frac{1}{k!}
\Vmatrix i_{m+1} &\ldots & i_{m+k}\\ 0 &\ldots &  0\endVmatrix-$$
$$-P_s\bmatrix i_1 &\ldots & i_m\\ j_1 &\ldots & j_m\endbmatrix
C^0_{s-(i_1+\dotsb+i_m+j_1+\dotsb+j_m)}\frac{1}{k!}
\Vmatrix i_{m+1} &\ldots &  i_{m+k}\\ 0 &\ldots & 0\endVmatrix=$$
$$=\cases 0,\quad\text{if $s\geqslant i_1+\dotsb+i_m+j_1+\dotsb+j_m,$}\\
P_s\bmatrix i_1 &\ldots & i_m\\ j_1 &\ldots & j_m\endbmatrix
\frac{1}{k!} \Vmatrix i_{m+1} &\ldots &  i_{m+k}\\ 0 &\ldots & 0\endVmatrix
, &\text{if $s<i_1+\dotsb+i_m+j_1+\dotsb+j_m$.$\square$}\endcases$$

\vskip 0.4cm
\noindent
{\bf 3. Equations for the Bacher-Akhiezer function}
\vskip 0.4cm

\noindent
Consider the KP hierarchy. This is a condition of compatibility of the
infinite system of the differential equations
$$\frac{\P\F}{\P x_n}=L_n\F,\ \tag 3.1$$
where $$L_n=\frac{\P^n}{\P x_1^n}+\sum^n_{i=2}B^i_n(x)
\frac{\P^{n-i}}{\P^{n-i} x_1},$$ and $\F$ is
a function of type
$$\F(x,k)=\text{exp}(\sum^\infty_{j=1}x_j k^j)(1+\sum^\infty_{i=1}
\X_i k^{-i}), $$
(here $k\in \Bbb C$ belong to some neighbourhood of $\infty$ and
$x=(x_1, x_2,...)$ --- is a  finite sequence).

Put us $$\P_i=\frac{\P}{\P x_i},\ \P=\P_1.$$ A direct
calculation gives

\vskip 0.4cm
\noindent
{\bf Lemma 3.1.} {\sl Conditions of compatibility of (1) are
$$B^t_s=-\sum^{t-1}_{i=1} C^i_s \P^i\X_{t-i}
-\sum^{t-1}_{j=2}B^j_s\sum^{t-j-1}_{i=0}C^i_{s-j} \P^i\X_{t-i-j}
, \tag 3.2$$
$$\P_n\X_i=\sum^{n+i-1}_{j=1} C^j_n \P^j\X_{i+n-j}+
\sum^{n}_{k=2}B^k_n \sum^{n-k}_{j=0}C^j_{n-k} \P^j\X_{i+n-j-k}
.\quad \square\tag 3.3$$ }

In this case $\F$ is called a {\sl Bacher-Akhiezer function}.

Consider now the function $$\ln\F(x,k)=\sum^\infty_{j=1} x_j k^j+
\sum^\infty_{j=1}\eta_j k^{-j},$$ where
$$\X_j=\sum^\infty_{n=1}\frac {1}{n!}\sum_{i_1+\dotsb+i_n=j}
\eta_{i_1}\dotsb \eta_{i_n}. $$

\vskip 0.4cm
\noindent
{\bf Lemma 3.2.} {\sl Let $2\leqslant t\leqslant s$.
Then $$B^t_s=-\sum_{n=1}^\infty \sum P_s\pmatrix i_1 &\ldots & i_n\\
j_1 &\ldots & j_n\endpmatrix \P^{j_1}\eta_{i_1}\dotsb \P^{i_n}\eta_{i_n},$$
where the second sum is taken by all matrices $\pmatrix i_1 &\ldots & i_n\\
j_1 &\ldots & j_n\endpmatrix$ such that $i_m\geqslant 1, j_m\geqslant 1$ and
$i_1+\dotsb +i_n+j_1+\dotsb +j_n=t$.}

\vskip 0.4cm
\noindent
Proof: An induction by $t$. For $t=2$ according to (3.2),
$B^2_s=-s\P\X_1=-P_s {1\choose 1}\P\eta_1$. Prove the lemma for
$t=N$ considering that it is proved for $t<N$. According to (3.2)
$$B^t_s=-\sum^{t-1}_{i=1} C^i_s \P^i
\Bigr(\sum^{\infty}_{n=1} \frac{1}{n!}\sum_{i_1+\dotsb +i_n=t-i}
\eta_{i_1}\dotsb\eta_{i_n}\Bigl)+$$
$$+\sum_{j=2}^{t-1}\Bigl(\sum^\infty_{n=1}\sum_{i_1+\dotsb +j_n=j} P_s
{\pmatrix i_1 &\ldots & i_n\\j_1 &\ldots & j_n\endpmatrix}
\P^{j_1}\eta_{i_1}\dotsb \P^{j_n}\eta_{i_n}\Bigr)\cdot$$
$$\cdot\Bigl(\sum_{i=0}^{t-j-1}C^i_{s-j}\P^j
\Bigl(\sum^\infty_{n=1}\frac {1}{n!}
\sum_{i_1+\dotsb +i_n=t-i-j}\eta_{i_1}\dotsb\eta_{i_n}\Bigr)\Bigr)=$$
$$=-\sum^\infty_{n=1}\sum\Bigl({1\over{n!}} C_s^{j_1+\dotsb +j_n}
\frac{(j_1+\dotsb +j_n)!}{j_1!\dotsb j_n!}-$$
$$-\sum^{n-1}_{q=1}P_s{\pmatrix i_1 &\ldots & i_q\\
j_1 &\ldots & j_q\endpmatrix}
\frac{1}{(n-q)!} C^{j_{q+1}\dotsb j_n}_{s-(i_1+\dotsb i_q+j_1+\dotsb
+j_q)}\frac {(j_{q+1}+\dotsb j_n)!}{j_{q+1}!\dotsb j_n!}\Bigr)
\P^{j_1}\eta_{i_1}\dotsb \P^{j_n}\eta_{i_n}=$$
$$=-\sum^\infty_{n=1}\sum P_s{\pmatrix i_1 &\ldots & i_n\\
j_1 &\ldots & j_n\endpmatrix}
\P^{j_1}\eta_{i_1}\dotsb \P^{j_n}\eta_{i_n},$$
where the second sums are taken by all matrices $\pmatrix i_1 &\ldots & i_n\\
j_1 &\ldots & j_n\endpmatrix$ such that $i_1+\dotsb +i_n+j_1+\dotsb +
j_n=t$, $i_m\geqslant 1, j_m\geqslant 0$.
According to lemma 2.1 it is possible to consider that in the last sum
$j_m>0$ for all $m$. $\square$

\vskip 0.4cm
\noindent
{\bf Lemma 3.3.} {\sl
$$\P_s\eta_r=\sum_{n=1}^\infty \sum P_s\pmatrix i_1 &\ldots & i_n\\
j_1 &\ldots & j_n\endpmatrix \P^{j_1}\eta_{i_1}\dotsb \P^{j_n}\eta_{i_n},$$
where the second sum is taken by all matrices $\pmatrix i_1 &\ldots & i_n\\
j_1 &\ldots & j_n\endpmatrix$ such that $i_m\geqslant 1, j_m\geqslant 1$ and
$i_1+\dotsb +i_n+j_1+\dotsb +j_n=r+s$.}

\vskip 0.4cm
\noindent
Proof: An induction by $r$. According to (3.3) and lemma 3.2 for
$r=1$ $$\P_s\eta_1=\P_s\X_1=
\sum^{\infty}_{j=1} C^j_s \P^j\X_{s+1-j}+
\sum^s_{k=2}B^k_s\sum^{s-k}_{j=0} C^j_{s-k}\P^j\X_{1+s-j-k}=$$
$$=\sum^s_{j=1} C^j_s\P^j\X_{s+1-j}-
\sum^{\infty}_{k=2}\Bigl(\sum^\infty_{n=1}\sum_{i_1+\dotsb +j_n=k} P_s
{\pmatrix i_1&\ldots & i_n\\j_1 & \ldots & j_n\endpmatrix}
\P^{j_1}\eta_{i_1}\dotsb \P^{j_n}\eta_{i_n}\Bigr)\cdot$$
$$\cdot\Bigl(\sum_{j=0}^{s-k}C^j_{s-k}\P^j\X_{1+s-j-k}\Bigr)=
\sum^\infty_{n=1}\sum P_s \pmatrix i_1 &\ldots & i_n\\
j_1 &\ldots & j_n\endpmatrix
\P^{j_1}\eta_{i_1}\dotsb \P^{j_n}\eta_{i_n},$$
where the second sum in the last formula is taken by all matrices
$\pmatrix i_1 &\ldots & i_n\\j_1 &\ldots & j_n\endpmatrix$ such that
$i_1+\dotsb +i_n+j_1+\dotsb +
j_n=s+1, \ i_m\geqslant 1, \ j_m\geqslant 0$.
According to lemma 1 in this sum it is sufficient consider
only matrices, where $j_m>0$ for all $m$.

Prove now the lemma for $r=N$, considering that it is proved for
$r<N$. According to (3.3)
$$\P_s\Bigl(\sum^\infty_{n=1}{\frac {1}{n!}}\sum_{i_1+\dotsb+i_n=r}
\eta_{i_1}\dotsb \eta_{i_n}\Bigr)=
\sum^{s+r-1}_{j=1}C^j_{s}\P^j\X_{s+r-j}+
\sum^\infty_{k=2}B^k_s\sum_{j=0}^{s-k} C^j_{s-k}\P^j\X_{r+s-j-k}.$$
Thus according to lemma 3.2, lemma 2.1 and inductive hypothesis,
$$\P_s\eta_r=\sum^{\infty}_{j=1}C^j_s\P^j
\Bigl(\sum^\infty_{n=1}{1\over n!}\sum_{i_1+\dotsb +i_n=s+r-j}
\eta_{i_1}\dotsb \eta_{i_n}\Bigr)+$$
$$+\sum^{\infty}_{k=2}\Bigl(\sum_{i_1+\dotsb +j_n=k} P_s
{\pmatrix i_1 &\ldots & i_n\\j_1 &\ldots & j_n\endpmatrix}
\P^{j_1}\eta_{i_1}\dotsb \P^{j_n}\eta_{i_n}\Bigr)\cdot$$
$$\cdot\sum_{j=0}^{s-k}C^j_{s-k}\P^j
\Bigl(\sum^\infty_{n=1}\frac {1}{n!}\sum_{i_1+\dotsb+i_n=r+s-j-k}
\eta_{i_1}\dotsb \eta_{i_n}\Bigr)
-\P_s\Bigl(\sum_{n=2}^{\infty}\frac{1}{n!}\sum_{i_1+\dotsb+i_n=r}
\eta_{i_1}\dotsb \eta_{i_n}\Bigr)=$$
$$=\sum^\infty_{n=1}\sum P_s
{\pmatrix i_1&\ldots & i_n\\j_1 &\ldots & j_n\endpmatrix}
\P^{j_1}\eta_{i_1}\dotsb \P^{j_n}\eta_{i_n},$$
where the second sum in the last formula is taken by all matrices
$\pmatrix i_1 &\ldots & i_n\\ j_1 &\ldots & j_n\endpmatrix$
such that $i_1+\dotsb +i_n+j_1+\dotsb +
j_n=s+r$, $i_m\geqslant 1, j_m\geqslant 1$. $\square$

\vskip 0.4cm
\noindent
{\bf 4. KP hierarchy}
\vskip 0.4cm

\noindent
According to [DKJM] the Bacher-Akhiezer function $\F$ is
$$\F(x,k)=\text{exp}(\sum x_j k^j)\frac{\tau(x_1-k^{-1},
x_2-\frac 12 k^{-2}, x_3-\frac 13 k^{-3},...)}{\tau(x_1,
x_2, x_3,...)}$$ for some function $\tau(x_1, x_2,...)$.
By analogy of [N1] this gives a possibility to
describe the KP hierarchy as an infinite system of
differential equations on $v(x,k)=-\ln \tau(x,k)$.
Really
$$\sum^\infty_{j=1}\eta_j k^{-j}=\ln\F(x,k)-\sum^\infty_{j=1}
x_j k^j=-v(x_1-k^{-1}, x_2-\frac12 k^{-2},...)+v(x)=$$
$$=\sum^\infty_{n=1}\sum_{i_1+\dotsb+i_n=j}\frac{(-1)^{n+1}}
{n!i_1\dotsc i_n}\P_{i_1}\dotsb \P_{i_n} v(x) k^{-j}.$$
Therefore
$$\eta_r=\sum^\infty_{n=1}\sum_{i_1+\dotsb+i_n=r}\frac{(-1)^{n+1}}
{n!i_1\dotsc i_n}\P_{i_1}\dotsb \P_{i_n} v. \tag 4.1$$

\vskip 0.4cm
\noindent
{\bf Theorem 4.1.} {\sl There exist universal rational
coefficients
$$R_r\pmatrix s_1 &\ldots & s_n\\ t_1 &\ldots & t_n\endpmatrix,
R_{ij}\pmatrix s_1 &\ldots & s_n\\t_1 &\ldots & t_n\endpmatrix$$
such that
$$\eta_r={1\over r}\P_rv+\sum^\infty_{n=1}\sum
R_r\pmatrix s_1 &\ldots & s_n\\t_1 &\ldots & t_n\endpmatrix
\P_{s_1}\P^{t_1}v\dotsb \P_{s_n}\P^{t_n}v , \tag 4.2$$
$$\P_i\P_j v=\sum^\infty_{n=1}
\sum R_{ij}\pmatrix s_1 &\ldots & s_n\\t_1 &\ldots & t_n\endpmatrix
\P_{s_1}\P^{t_1}v\dotsb \P_{s_n}\P^{t_n}v, \tag 4.3$$
where the second sums are taken by all matrices
$\pmatrix s_1 &\ldots & s_n\\ t_1 &\ldots & t_n\endpmatrix$ such
that $s_m, t_m\geqslant 1$, and the sum $s_1+\dotsb +s_n+t_1+
\dotsb +t_n$ is equal $r$ for (4.2) and $i+j$ for (4.3).}

\vskip 0.4cm
\noindent
Proof: An induction by $k$  and $i+j$. For $i+j=2$ the theorem
is obviously. For $r=1$ it follows from (4).
Prove the theorem for $i+j=N$ and $r=N-1$, considering that it is
proved for $i+j<N$ and $r<N-1$.
Later we consider that $s_m, t_m\geqslant 1$ and
$\SI_n=s_1+\dotsb +s_n+t_1+\dotsb +t_n$.
Then according to (4.1) and (4.3)
$$\eta_r={1\over r}\P_rv+\sum^\infty_{n=2}\sum_{s_1+\dotsb+s_n=r}
\frac{(-1)^{n+1}}{n! s_1\dotsb s_n}\P_{s_1}\dotsb \P_{s_n}v(x)=$$
$$={1\over r}\P_rv+\sum^\infty_{n=1}\sum_{\SI_n=r}
R_r\pmatrix s_1 &\ldots & s_n\\t_1 &\ldots & t_n\endpmatrix
\P_{s_1}\P^{t_1}v\dotsb \P_{s_n}\P^{t_n}v.$$
Thus according to (4.2), (4.3) and lemma 3.3
$${1\over j}\P_i\P_j v=\P_i\eta_j-\P_i
\Bigl(\sum^\infty_{n=1}\sum_{\SI_n=j}
R_j\pmatrix s_1 &\ldots & s_n\\t_1 &\ldots & t_n\endpmatrix
\P_{s_1}\P^{t_1}v\dotsb \P_{s_n}\P^{t_n}v\Bigr)=$$
$$=\sum^\infty_{n=1}\sum_{\SI_n=i+j}
P_i\pmatrix s_1 &\ldots & s_n\\t_1 &\ldots & t_n\endpmatrix
\P^{t_1}\eta_{s_1}\dotsb \P^{t_n}\eta_{s_n}-$$
$$-\P_i\Bigl(\sum^\infty_{n=1}\sum_{\SI_n=j}
R_j\pmatrix s_1 &\ldots &s_n\\t_1 &\ldots & t_n\endpmatrix
\P^{t_1}\P_{s_1}v\dotsb \P^{t_n}\P_{s_n}v\Bigr)=$$
$$=\sum^\infty_{n=1}\sum_{s_1+\dotsb+ s_n+n=i+j}
P_i\pmatrix s_1 &\ldots & s_n\\1 &\ldots & 1\endpmatrix
\P(\frac{1}{s_1}\P_{s_1}v)\dotsb \P(\frac{1}{s_n}\P_{s_n}v)+$$
$$+\sum^\infty_{n=1}\sum_{\SI_n=i+j,t_1+\dotsb+ t_n>n}
R_{ij}\pmatrix s_1 &\ldots & s_n\\t_1 &\ldots & t_n\endpmatrix
\P_{s_1}\P^{t_1}v\dotsb \P_{s_n}\P^{t_n}v. \square $$

\noindent
{\bf Remark 4.1.} The system (4.3) was at first deduced in
[DN]. The set of its solution bijectively
correspond to the set of solutions of (3.1). Up to a constant
formal solution of (4.3) is defined by an infinite set of
functions of one variable
$f_i(x_1)=\P_i v|_{x_2=x_3=\dotsb=0}\quad (i=1,2,...)$.

\noindent
{\bf Remark 4.2.} The algorithm described in the proof ot theorem 4.1
gives an algorithm for calculation of all rational constants
$R_{ij}\pmatrix s_1 &\ldots & s_m\\t_1 &\ldots & t_m\endpmatrix$.
The first equations of hierarchy (4.3) are:
$$\P_2^2 v=\frac 43\P_3\P v-\frac 13\P^4 v+2(\P^2 v)^2, \tag 4.4$$
$$\P_3\P_2 v=\frac 32\P_4\P v-\frac 32\P_2\P^3 v+3\P_2\P v
\P^2 v,$$ $$\P_3^2 v=\frac 95\P_5\P v-\P_3\P^3 v+\frac 15
\P^6 v+3\P_3\P v \P^2 v+\frac 94(\P_2\P v)^2-3\P^4 v\P^2 v-
\frac 94(\P^3 v)^2+3(\P^2 v)^3.$$
The equation (4.4) is KP equation twice integrated over $x_1$.

\vskip 0.4cm
\noindent
{\bf Theorem 4.2.} {\sl If $\sum^m_{i=1}(t_i+1)\equiv 1
(\text{mod 2})$, then $R_{ij}\pmatrix s_1 &\ldots & s_m\\t_1
&\ldots & t_m\endpmatrix =0$.}

\vskip 0.4cm
\noindent
Proof: The equations of KP--hierarchy are equivalents of equations on
function $\tau(x)$ [DKJM]. All these equations can be written
simply by means of the "bilinear Hirota operators". We recall
the definition of them. If $f(x)$ is a function of one variable, then
for any polynomial (or power series) $Q$ the action of the Hirota
operator $Q(D_x)f(x)\cdot f(x)$ is defined by
$$Q(D_x)f(x)\cdot f(x)=Q(\P_y)[f(x+y)f(x-y)]_{y=0}.$$
For functions of several variables the definition is similar.
The generating function for the equations of the KP hierarchy
has the form
$$\sum^\infty_{j=0} p_j(-2y)p_{j+1}(\widetilde D)\exp\big(
\sum^\infty_{i=1} y_i D_i\big)\tau\cdot\tau=0, \tag 4.5$$
where $y=(y_1, y_2,\dots,)$ are auxiliary independent variables,
$\widetilde D= (D_1, \allowbreak 2^{-1} D_2, 3^{-1} D_3,\dots)$,
$D_j$ is the Hirota operator in the  variable $x_j$ and $p_j$
are the Schur polynomials defined from the following expansion:
$$\exp\big(\sum^\infty_{j=1} x_j k^j\big)=\sum^\infty_{j=0}
k^j p_j(x_1,\dots, x_j).$$

All monomials of odd degree give trivial Hirota operators.
Therefore if $\tau(x)$ is a solution of (4.5), then $\widetilde \tau
(x)=\tau(-x)$ is also solution of (4.5). Moreover, according to [DN]
a function $\tau$ is a solution of the system (4.5) if and only if
$v=-\ln(\tau)$ is a solution of the system (4.3). Thus,
$v(x)$ is a formal solution of the system (4.3), if and only if
$\widetilde v(x)=v(-x)$ is a formal solution of the system
(4.3). This is equivalent of the affirmation of theorem 4.2.
$\square$.

\vskip 0.6cm
\noindent
{\bf 5. Gelfand--Dikii hierarchy}
\vskip 0.4cm

\noindent
According to [S] the set of solution of $n$--Gelfand--Dikii hierarchy
bijectively correspond to the set of nondepending from $x_n$
solutions of KP hierarchy. In this case according to theorem 4.1
$$0=\P_m\P_n v=\frac{mn}{m+n-1}\P_{n+m-1}\P v+\sum_{m=1}^\infty\sum
R_{mn}\pmatrix s_1&\ldots&s_m\\t_1&\ldots&t_m\endpmatrix
\P_{s_1}\P^{t_1}v\dotsb\P_{s_m}\P^{t_m}v,$$ where
$1\leqslant s_j\leqslant n+m-2,\ t_j\geqslant 1$. This gives recurrence
formulas expressing $\P_k\P v$ for $k>n$ via $\P_r\P v$ for $r<n$.
Thus we have relations
$$\P\P_{n+r} v=\sum_{m=1}^\infty\sum N_{1(n+1)}^m
\pmatrix s_1&\ldots&s_m\\t_1&\ldots&t_m\endpmatrix
\P_{s_1}\P^{t_1}v\dotsb\P_{s_m}\P^{t_m}v, \tag 5.1$$ where
$t_j\geqslant 1, s_j<n,\quad \sum^m_{j=1}(s_j+t_j)=n+r+1$.

\noindent
{\bf Example.} For $n=2$ the system (5.1) passes to $K\P V$
hierarchy.

Compering the systems (5.1) and (4.3) we find the system
$$\P_i\P_j v=\sum_{m=1}^\infty\sum
N_{i_j}^m\pmatrix s_1&\ldots&s_m\\t_1&\ldots&t_m\endpmatrix
\P_{s_1}\P^{t_1}v\dotsb\P_{s_m}\P^{t_m}v,\tag 5.2$$ where
$i,j\geqslant 1,\
1\leqslant s_\A\leqslant n-1,\ t_\A\geqslant 1$,
$\sum_{i=1}^m(s_\A+t_\A)=i+j$ and $N_{ij}^m
\pmatrix s_1&\ldots&s_m\\t_1&\ldots&t_m\endpmatrix$ are
some universal rational coefficients.

\noindent
{\bf Examples.}

1. For $n=3$ the first equation from the system (5.2) is
the Boussinesq equation
$$\P^2_2 v=-\frac 13\P^4 v+2(\P^2 v)^2.$$
2. For $n=4$ the first equations of  the system (5.2) are
$$\P_2^2 v=\frac 43\P_3\P v-\frac 13\P^4 v+2(\P^2 v)^2,$$
$$\P_3\P_2 v=-\frac 32\P_2\P^3 v+3\P_2\P v \P^2 v,$$
$$\P_3^2 v=-\frac 14\P_3\P^3 v+\frac 18\P^6 v+\frac 98
(\P_2\P v)^2-\frac 98(\P^3 v)^2-\frac 94\P^4 v\P^2 v+ 3(\P^2 v)^3.
\square $$

\vskip 0.4cm
\noindent
{\bf Theorem 5.1.} {\sl The Gelfand--Dikii hierarchy is equivalent to
a system of differential equations in a form
$$\P_{i_1}\dotsb\P_{i_k}v=\sum_{m=1} \sum_{i_1\dotsb i_k}
N_{i_1\dotsb i_k}^m\pmatrix s_1&\ldots&s_m\\t_1&\ldots&t_m\endpmatrix
\P_{s_1}\P^{t_1}v\dotsb\P_{s_m}\P^{t_m}v, \tag 5.3$$ where
$k\geqslant 2, t_j\geqslant 1, s_j<n,\ \sum_{j=1}^ki_j=\sum^m_{j=1}(s_j+t_j),
\ \sum^m_{j=1}t_j\geqslant m+k-2$ É $k+m+\sum^m_{j=1}t_j\equiv
0\ (\text{mod} 2)$.}

\vskip 0.4cm
\noindent
Proof: For $k=2$ the equations (5.3) coincide with the equations
(5.2). For $k>2$ the equations (5.3) are received from equations for
$k-1$ by differentiation by $\P_{i_k}$ and replacing
$\P_i\P_{i_k}v$ by (5.2). This gives monomials
$\P_{s_1}\P^{t_1}\dotsb\P_{s_m}\P^{t_m}v$, where the $\sum^k_{j=1}i_j=
\sum^m_{j=1}(s_j+t_j)$ É $\sum^m_{j=1}t_j\geqslant m+k-2$. The
condition $k+m+\sum^m_{j=1}t_j\equiv 0(\text{mod} 2)$ follows from
theorem 4.2. $\square$

\noindent
{\bf Remark 5.1.} Structure of the system (5.3) such that its formal
solutions are defined up to constant by arbitrary set of $n-1$
series from one variable  $f_i(x_1)=
\P_iv\vert_{x_2=x_3=\dotsb =0}\quad (i=1,...,n-1)$.

\noindent
{\bf Remark 5.2.} Coefficients
$N_{i_1\dotsb i_k}^m\pmatrix s_1&\ldots&s_m\\t_1&\ldots&t_m\endpmatrix$
are rational cons\-tants. The construc\-tions, described in the proofs
of theorems  4.1 É 5.1, give recurrent formulas for its calculation.

\vskip 0.4cm
\noindent
{\bf 6. Witten solution of the Gelfand--Dikii hierarchy}
\vskip 0.4cm

\noindent
Follow by Witten [W] let us consider the space $M_{g,s}$ of Riemann
surfaces of genus $g$ with $s$ punctures. Correspond to any puncture
a pair $(k_i, m_i)$, where $1\leqslant k_i<n,\quad m_i
\geqslant 0$. Witten [W] connects with the set $\{(k_i, m_i)
\vert i=1,...,s\}$ a number (correlator) $<\prod\limits_{k,m}
\tau^{d_{k,m}}_{k,m}>_g$, where $d_{k,m}$ is the number of pairs
$(k_i, m_i)$, that equal to $(k, m)$. The number
$<\prod\limits_{k,m} t^{d_{k,m}}_{k,m}>_g$ is equals to the value of
some class of cohomology on a compactification of a space of
$n$--spin bundles over $P\in M_{g,s}$ [W]. Put us
$$F^g(t_{1,0}, t_{1,1},...)=\sum_{d_{k,m}}<\prod\limits_{k,m}
\tau^{d_{k,m}}_{k,m}>_g\prod\limits_{k,m}
\frac{t^{d_{k,m}}_{k,m}}{d_{k,m}!}.$$
According to the Witten conjecture the serias $F=\sum\limits^\infty_{g=0}
F^g$ after the change $t_{k,m} \mapsto -(mn+k)x_{mn+k}$ pass to
a formal solution $v$ of the system (5.3) satisfying the equation
$$\P v=\frac 12\sum_{i+j=n}ij x_ix_j+\sum_{i=1}^\infty
(i+n)x_{i+n}\P_i v, \tag 6.1$$
$$0=v(0)=\P_i v(0)\quad (i=1,2,...).$$
The single such solution $W$ we call the Witten solution of the
Gelfand--Dikii hierarchy.

\vskip 0.4cm
\noindent
{\bf Theorem 6.1.} {\sl The Witten solution of the Gelfand--Dikii
hierarchy is $W=\sum\limits^\infty_{g=0}W^g$, where
$$W^g(x_1, x_2,...)=$$ $$=\sum^\infty_{k=2}\sum_{i_1+
\dotsb+i_k=(n+1)(2g-2+k)}\frac{(n-1)^{2g-2+k}}{k!}N^{2g-2+k}_{i_1
\dotsb i_k} \pmatrix n-1&\ldots&n-1\\2&\ldots&2\endpmatrix
x_{i_1}\dotsb x_{i_k}.$$ Moreover the function $W^g$ is a
quasihomogeous series of degree $(1-g)(2+\frac 2n)$ by $x_i$
of degrees $1+\frac 1n-\frac in$.}

\vskip 0.4cm
\noindent
Proof: Compatibility of the equations (5.3) É (6.1) follows from
[AM], where this solution is represented in a form of matrix
integral. According to (6.1) $\P\P_i W\vert_{x_2=x_3=\dotsb=0}=
\delta_{n-1,i}(n-1)x_1$. These conditions and the equations (5.3)
uniquely determine all functions
$f_{i_1\dotsb i_k}(x_1)=\P_{i_1}\dotsb\P_{i_k} W\vert_{x_2=x_3=
\dotsb=0}$. According to theorem 5.1 if $f_{i_1\dotsb i_k}(0)\ne 0$,
that $(i_1+\dotsb+i_k)\equiv 0\ (\text{mod} (n+1))$ and
$$f_{i_1\dotsb i_k}(0)=(n-1)^m N^m_{i_1\dotsb i_k}
\pmatrix n-1&\ldots&n-1\\2&\ldots&2\endpmatrix,$$ where
$m=\frac{i_1+\dotsb+i_k}{n+1}$. From this theorem follow also that
$m\geqslant k-2$ É $k+m\equiv 0\ (\text{mod} 2)$. Therefore
$m=2g+k-2$, where $g\geqslant 0$ is a natural number. Thus,
$W=\sum W^g$, where
$$W^g(x_1, x_2,...)=\sum^\infty_{k=2}\sum_{i_1+
\dotsb+i_k=(n+1)(2g+k-2)}\frac{1}{k!}f_{i_1\dotsb i_k}(0)x_{i_1}
\dotsb x_{i_k}=$$
$$=\sum^\infty_{k=2}\sum_{i_1+\dotsb+i_k=(n+1)(2g+k-2)}\frac{1}{k!}
(n-1)^{2g-2+k}N^{2g-2+k}_{i_1\dotsb i_k} \pmatrix n-1&\ldots&n-1\\
2&\ldots&2\endpmatrix x_{i_1}\dotsb x_{i_k}.$$
The quasihomogeneity of the series $W^g$ follows from $i_1+\dotsb
+i_k=(n+1)(2g-2+k)$. $\square$

\vskip 0.4cm
\noindent
{\bf Corollary 6.1.} {\sl The Witten solution of the
Gelfand--Dikii hierarchy has a representation by the sum of
quasihomogeneous series $W^g$ of the same degrees that $F^g$.}

\vskip 0.4cm
\noindent
Proof: According to [W] $F^g$ is a quasihomogeneous series of
degree $(1-g)(2+\frac 2n)$ by $t_{k,m}$ of degree
$1+\frac 1n-k-\frac mn$. $\square$

\vskip 0.4cm
\noindent
{\bf Theorem 6.2.} {\sl The functions $W$ and $W^0$ coincide on
the set $L_0=(x_1, x_2,\dots, x_{n-1},
\allowbreak 0,0,\dots)$.}

\vskip 0.4cm
\noindent
Proof: According to (6.1) $\P^r\P_\ell W=0$ on the set
$L_0$ if $\ell<n-1,\ r>1$, or $\ell= n-1,\ r>2$.
Besides according to (5.1)
$$\P\P_{n+\ell} W=\sum^\infty_{m=1}\sum  N^m_{1(n+\ell)}
\pmatrix s_1 &\ldots & s_m\\t_1 &\ldots & t_m\endpmatrix
\P_{s_1}\P^{t_1}W\dotsb \P_{s_m}\P^{t_m}W,$$
where $\sum\limits_{i=1}^m(s_i+t_i)=n+\ell+1$. Thus, if
$\ell<n$ and $\P\P_{n+\ell}W \ne 0$,
that all numbers $t_j$ are less than 3 and
among them is not  two number more than 1. But if among the numbers
$1\leqslant t_1,\dots,t_m\leqslant 2$ there is exactly one $t_j=2$
then $N^m_{1(n+\ell)}\pmatrix s_1 &\ldots & s_m\\t_1 &\ldots &
t_m\endpmatrix= 0$ by theorem 5.1. Thus,
$$\P\P_{n+\ell} W=\sum^\infty_{m=1}\sum
N^m_{1(n+\ell)}\pmatrix s_1 &\ldots & s_m\\1 &\ldots &
1\endpmatrix \P_{s_1}\P W\dotsb \P_{s_m}\P W.$$
Moreover, according (6.1) $\P_\ell W=\P\P_{n+\ell}W$ and
$\P_s\P W=s(n-s)x_{n-s}$ on the set $L_0$ . Thus from the
equality $\sum\limits^m_{i=1}(s_i+t_i)=n+\ell+1$
follows that $\P_\ell W\vert_{L_0}$ is a quasihomogeneous
polynomial of degree $2+\frac 2n-(1+\frac 1n-\frac\ell n)=
1+\frac 1n+\frac\ell n$. Thus, $W\vert_{L_0}$ is a quasihomogeneous
polynomial of degree $2+\frac 2n$ É $W\vert_{L_0}=W^0\vert_{L_0}.
\square$

\vskip 0.4cm
\noindent
{\bf Corollary 6.2.} {\sl $F(x_1, x_2, \dots, x_{n-1}, 0, 0, \dots)=
W(x_1, \frac{x_2}{2}, \frac{x_3}{3}, \dots, \frac{x_{n-1}}{n-1},
0, 0,\dots)$.}

\vskip 0.4cm
\noindent
Proof: According to [W] $F^0(x_1, x_2, \dots)=
W^0(x_1, \frac{x_2}{2}, \frac{x_3}{3}, \dots)$ and
$F(x_1, x_2, \dots, x_{n-1}, 0,\allowbreak 0,\dots)=
F^0(x_1, x_2, \dots, x_{n-1}, 0,0,\dots)$. Thus
the theorem 6.2 imply corollary 6.2.  $\square$

\vskip 0.4cm
\noindent
{\bf Remark 6.1.} According to [DVV, Kr] the function
$W^0\vert_{L_o}$ is the potential of Frobenius structure on the
space of versal deformations of the singularity $A_n$. By theorem
6.2 we have $W^0\vert_{L_0}=W\vert_{L_0}$. A simple
algorithm of calculation of this function describes in [N2].

\vskip 0.4cm
\noindent
{\bf References}

[AM] M.Adler, P.van Moerbeke, A matrix integral solution to
two--dimensional $W_p$--gravity. Commun. Math. Phys. 147 (1992),
25-56.

[DKJM]  E.Date, M.Kashiwara, M.Jimbo, T.Miwa, Transformation groups
for so\-li\-ton equa\-tion, Proceedings
of RIMS Sym\-po\-sium on Non-Linear Integrable System.
Sin\-gapo\-re: World Science Publ. Co., 1983, 39-119.

[DN]  B.A.Dubrovin, S.M.Natanzon, Real theta-function solutions
of the Ka\-domt\-sev--Petvi\-ashvili equation. Math. USSR Irvestiya,
32:2 (1989), 269-288.

[DVV]  R.Dijkgraaf, E.Verlinde, H.Verlinde, Topological strings
in $d<1$. Nucl. Phys., B 352 (1991), 59.

[GD] I.M.Gelfand, L.A.Dikii, The resolvent and Hamiltonian
systems, Funct. Anal. Appl. 11 (2) (1977), 93-105.

[Ko] M.Kontsevich, Intersection theory on the moduli space of
curves and the matrix airy function. Commun. Math. Phys. 147
(1992), 1-23.

[Kr]  I.Krichever, The dispersionless Lax equations and topological
minimal models. Com\-mun. Math. Phys. 143 (1992), 415-429.

[N1]  S.M.Natanzon, Real nonsingular finite zone solutions of
solution equations. Amer. Math. Soc. Transl. (2). V.170 (1995),
153-183

[N2] S.M.Natanzon, Formulas for $A_n$ and $B_n$--solutions of WDVV
equations, \linebreak hep-th/9904103.

[S] M.Sato, Soliton equations and universal Grassmann manifold.
Math. Lect. Notes Ser., Vol.18, Sophia University, Tokyo (1984).

[W]  E.Witten, Algebraic geometry associated with matrix models
of two di\-men\-sio\-nal gravity, Topological models in modern mathematics
(Stony Brook, NY, 1991), Publish or Perish, Houston, TX 1993, 235-269.

\end